# Major Maintenance Schedule Optimization for Electric Multiple Unit Considering Passenger Transport Demand

Jianping Wu[1], Boliang Lin[1*],

[1] School of Traffic and Transportation, Beijing Jiaotong University, Beijing 100044, China

15114205@bjtu.edu.cn, bllin@bjtu.edu.cn

**Abstract**

It is an important objective pursued in a railway agency or company to reduce the major maintenance costs of electric multiple unit (EMU). The EMU major maintenance schedule decides when to undergo major maintenance or undertake transportation task for train-set, based on practical requirements, such as passenger transport demand, workshop inspection capacity, and maintenance requirements. Experienced railway practitioners can generally produce a feasible major maintenance schedule; however, this manual process is time-consuming, and an optimal solution is not guaranteed. This research constructs a time-space network that can display the train-set status transformation process between available and major maintenance status. On this basis, a 0-1 integer programming model is developed to reduce the major maintenance costs with consideration of all necessary regulations and practical constraints. Compared with the manual process, the genetic algorithm with simulated annealing survival mechanism is also developed to improve solution quality and efficiency. It can reduce the complexity of the algorithm substantially by excluding infeasible solutions when constructing the model.

**Introduction**

According to the statistical bulletin of National Railway Administration of the People's Republic of China, the railway passenger traffic volume reached 295 million during the Spring Festival Travel Rush in 2015 (from Feb.4 to Mar.15, 40 days in total), among which people travelling by high-speed rail accounts for 41.4%. China Railway High-speed (CRH) operated by China Railway Corporation (CR) is more and more popular for its convenience and comfortableness. How to ensure enough rolling stock to fulfill heavy transportation tasks at rush hours (Spring Festival Travel Rush in particular) has long been the difficult problem for CR. The high purchase cost and complicated maintenance system of electric multiple unit (EMU), which is the unique vehicle running on CRH, makes the problem even worse. CRH runs different EMU train-sets, the designs for which are imported from other nations and designated CRH1 through CRH5 and CRH380A(L), CRH380B(L), CRH380C(L), CRH380D(L) and so on.[1]

In CR system, regular preventive maintenance mechanism is implemented for CRH series train-sets. There are five levels of maintenance, e.g., first-level (daily inspection), second-level (special inspection), third-level (bogie inspection), fourth-level (system decomposition inspection) and fifth-level (general inspection). Maintenances at the first and second level are both operational inspection with short inspection cycle and inspection time. Typically, the maintenance plans of these two levels are co-optimization with the rolling stock operational plan respectively. The rest belong to major maintenance. Major maintenance are scheduled in advance for each train-set during tactical maintenance planning, which are drawn up once a year. That is because it requires comparatively longer inspection cycle, longer inspection time and the limited inspection capacity. In general, the major maintenance schedule starts from September or October, and ends before the second Spring Festival Travel Rush coming. The whole planning horizon lasts for about 18 months or so. Therefore, the planning horizon is longer than the interval of planning, and some train-sets belong to this planning horizon will go to workshop for inspecting during the next schedule. For example, major maintenance schedule starts from Oct. 18th every year (e.g. 2015), and some train-sets will go to workshop on Nov. 20th, 2016. So, the maintenance date of these train-sets will be adjusted according to actual situation or not, and then, as the known conditions input the next schedule.

Because of the major maintenance is time-consuming, it will lead to the lack of train-sets for undertaking transportation tasks when the process of major maintenance be undergone at rush hours. However, major maintenance has long cycle, its arrangement is flexible. Therefore, this paper focuses on the optimization of electric multiple unit major maintenance schedule (EMUMMS), trying to avoid rush hours (especially the Spring Festival Travel Rush), and takes the cost of major maintenance into account to prevent frequent maintenance, which is of great significance in supply enough train-set for passenger transport peak demand. In the meantime, it also provides a decision support for capacity configuration of workshop.

**Literature review**

For the maintenance system, Stuchly et al.[2] analyzed the maintenance management systems, and developed a maintenance management system based on reliability-centered maintenance which could help engineers to schedule the maintenance scheme reasonably according to the technical conditions reflected by the EMU real-time monitoring data. Shimada[3] introduced a new maintenance system of accident prevention based on the maintenance technology situation of Japanese Shinkansen, thereby reducing the maintenance costs and improving the utilization efficiency of EMU. Cheng and Tsao[4] proposed a selection strategy of EMU maintenance based on the characteristics of the preventive and corrective maintenances. The estimation method for the spare parts' quantities and replacement intervals of specific components of EMU were also provided.

For the optimization of EMUs maintenance schedule, the existing literature is more concerned about the first- and the second-level maintenances. Sriskandarajah et al.[5] optimized the EMU overhaul plan using the genetic algorithm. The genetic algorithm improves the quality of plan and reduces the costs of operation. Maróti and Kroon[6,7] developed a multicommodity flow model for preventive maintenance routing. Given the EMU require maintenance in the forthcoming one to three days, the operation schedule must be adjusted so that these urgent units reach the maintenance facility in time. Alfieri et al.[8] proposed a multicommodity flow model for efficient rolling stock circulation on a single line of the Dutch railway network. Their objective is to minimize the distance run by train units of various types. Short-term maintenance requirements are not considered in their formulations. Rezvanizaniani et al.[9] discussed the implementation of Reliability Centered Maintenance to make rolling stock maintenance of the Raja Passenger Train Corporation more cost effective by reducing erroneous maintenance and unnecessary maintenance. Tsuji et al.[10] analyzed the influence factors of EMU operation and maintenance problems combined with the Japanese Shinkansen. And they developed a novel approach based on ant colony optimization to solve the problems. Wang et al.[11] proposed a 0-1 integer programming model to study EMU maintenance, but the model cannot be applied in major maintenance practice. Giacco et al.[12] provided a mixed-integer linear-programming formulation for integrating short-term maintenance planning in a network-wide railway rolling stock circulation problem, and the optimization objective is to minimize the cost with service pairings, empty runs, and short-term maintenance tasks. Lai et al.[13] studied an exact optimization model to improve the efficiency in rolling stock usage with consideration of practical operation constraints and multi-level maintenance. Compared to the manual process, a hybrid heuristic process is developed to improve solution quality and efficiency.

For the locomotive, Ziarati et al.[14] studied the locomotive operation by constructing a large integer programming model, they focused on the influence of maintenance for operation. Lingaya et al.[15] described a model for operational management of locomotive hauled railway cars. They sought for a maximum expected profit schedule that satisfied various constraints, among them maintenance requirements were also included. Wang et al.[16] developed a 0-1 integer programming model to study locomotive operation and maintenance, and applied genetic algorithm to solve the model. But the model ignored the maintenance capacity.

In the optimization of planes maintenance schedule research, Moudania and Félix[17] studied the problems of assigning planes to flights and of fleet maintenance operations scheduling. And they also proposed a dynamic approach to cope with the fleet assignment problem and a heuristic technique to solve the embedded maintenance scheduling problem. Budai et al.[18] presented a mathematical formulation for the long-term planning of railway maintenance works. The objective is to minimize the time required for maintenance, expressed as a cost function. Heuristic algorithms computed nearly optimal solutions by combining maintenance activities on each track. Mehmet and Bilge[19] developed integer linear programming model by modifying the connection network representation which provided feasible maintenance routes for each aircraft in the fleet over a weekly planning horizon with the objective of maximizing utilization of the total remaining flying time of fleet. The proposed model is solved by using branch-and-bound under different priority settings for variables to branch on.

In addition, Grigoriev et al.[20] studied the problem of scheduling maintenance services, and the objective was to find a cyclic maintenance schedule of a given length that minimized total operation costs. A branch and price algorithm was applied to solve the model.

Compared to the existing researches, this study makes the some new contributions. The goal of this paper is to explore the optimization of EMUMMS to meet passenger flow peak demand. Optimize electric multiple unit major maintenance schedule by constructing a time-space network. And a 0-1 integer programming model is proposed to reduce the costs of major maintenance with consideration of all necessary regulations and practical constraints, especially the passenger transport peak demand. An effective and intelligent genetic algorithm with simulated annealing survival mechanism is also developed for a general large-scale empirical case.

**EMU major maintenance problem at CR**

The EMUMMS is a plan that decides the train-set when to undergo major maintenance and when to undertake transportation task. The planner should not only ensure supply enough available train-set for passenger transport peak demand, but also reduce the maintenance cost as much as possible. In practice, steps of making an EMUMMS are as follows. First, collect the real-time data of all train-sets, such as the total operation mileage and days since last major maintenance, the total operation mileage and days, the level of last maintenance, start date of scheduled, average daily operation mileage of different type's train-sets and so on. Second, estimate the start date of major maintenance for all train-sets based on the statistics data, and maintenance requirements. At last, with the all necessary regulations and practical constraints, determine the exact start date of major maintenance for all train-sets. These are the outcome we needed. In CR system, an EMU usually consists of eight or sixteen units of power and non-power rolling stock permanently connected together, we mark it 8-EMU or 16-EMU. Therefore, each unit within a train-set operates daily tasks and undergoes inspections together. Depending on the demand of train trips, a train can be formed by one to two 8-EMU or one 16-EMU.

**Demand**

With the development of China's economy, the improvement of people's living standard as well as a huge change of consumption concept, passenger trip has formed a trend of long travel rush e.g. Spring Festival, Summer Holiday and National Day, coupled with short travel rush, such as New Year's Day, Tomb-sweeping Day, Labor Day, Dragon Boat Festival and Mid-autumn Day. Aimed at these travel rushes, China Railway Corporation has made passenger train graph based on relevant documents and files, and every railway bureaus have drawn up corresponding train-set operational plans. The EMUMMS is to ensure supply enough well-conditioned train-sets to implement the operational plan (particularly in travel rush).

**Maintenance/Inspection**

Table 1 shows the major maintenance requirement for the CRH2 and CRH380A(L) in terms of cumulative operating mileage and cumulative operating days.

**Table 1.** Maintenance regulations for CRH2 and CRH380A in CR

| Type | Third-level | Fourth-level | Fifth-level |
| --- | --- | --- | --- |
| CRH2A |  |  |  |
| CRH2B |  |  |  |
| CRH2C-1 | 0.6 million$_{-50,000}^{+20,000}$ km | 1.2 million$_{-100,000}^{+50,000}$ km | 2.4 million ±100,000 km |
| CRH2C-2 | or 1.5 years | or 3 years | or 6 years |
| CRH380A |  |  |  |
| CRH380AL |  |  |  |

The maintenance requirement sets the limit on how much distance and time a train-set can operate before the next mandatory maintenance. For example, the train-set, CRH2A, must go through the third-level maintenance after operating for 0.6 million km or for 18 months before it can be assigned to the next operation utilization path with actual train trips. Furthermore, the major maintenance requirement for train-set has a floating range. For example, cumulative operating mileage floating range of the third-level maintenance for CRH2A is [-50,000km, +20,000km], i.e., CRH2A will be undergone the procedures of third-level maintenance when its accumulated operation mileage is greater than 550.000 km and less than 620,000 km from the last major maintenance. And, the interval between two adjacent major maintenances should be less than the third-level requirements.

According to Table 1, the fourth-, fifth-level cycle is twice and four times, respectively, than the third-level cycle. Moreover, the maintenance procedures of higher-level maintenance include all maintenance procedures in lower-level maintenance. Thus, after one class of maintenance process, all cumulative operating mileage and days associated with that class and the corresponding lower-level classes of maintenance are set to zero. For example, all associated cumulative operating mileage and days for fourth-level and third-level should be set to zero after fourth-level. Therefore, the major maintenance requirements equals to the third-level requirements, and the train-set undergoes two adjacent major maintenances at different levels accordingly. On the basis of this rule, the major maintenance cyclic graph could be drawn, as shown in Fig. 1. Based on the sequence of major maintenance cyclic graph and the major maintenance requirement, the train-set accepts the different levels of major maintenances in turn until they are scrapped.

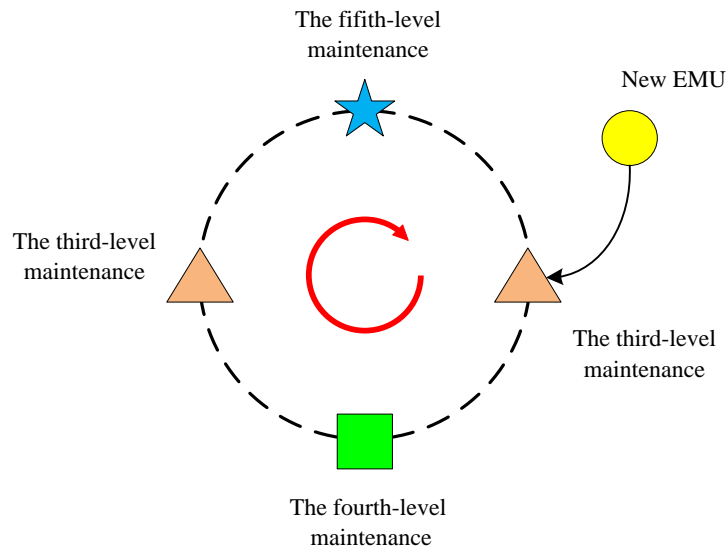

Fig. 1. The major maintenance cyclic graph

According to *User Manual* [21], the major maintenance requirements depends on the cumulative operating mileage mainly and the cumulative operating days as a supplement, and whichever comes first. According to the simple mathematical calculation based on the data in Table 1 and Table 2, we can conclude that the cumulative operating mileage for any train-set will finish firstly, such as CRH2A, (0.6 million km)/ (1500 km) =400 days < 1.5 years. Therefore, we could only take into account the cumulative operating mileage requirement only and ignore the cumulative operating days in this paper.

**Table 2.** Average daily operation mileage of different type's train-sets

| Type | CRH2A | CRH2B | CRH2C-1 | CRH2C-2 | CRH380A | CRH380AL |
|---|---|---|---|---|---|---|
| Average daily mileage | 1,500 km | 1,500 km | 1,600 km | 1,800 km | 1,900 km | 1,900 km |

When the train-set needed major maintenance, it would be detained and sent to workshop. Since the major maintenance is time-consuming, it will result in shortage of well-conditioned train-set when the major maintenance be undergone during the peak period of passenger flow (such as Chinese "Spring Festival travel rush") with the fleet-size constraints. But the floating range ensures the start time of maintenance is flexible. In order to avoid the peak period of passenger flow, some train-sets may choose to undergo the major maintenance at the rear of floating range, and others may conduct maintenance at the head of floating range. However, too much delay will influence operation safety, and too much maintenance in advance will increase the frequency. In addition, as the major maintenance is costly, frequent maintenance will lead to a great waste of resources. Therefore, under the premise of safety operation, the start time of major maintenance should be postponed. Aside from these requirements on rolling stock, every depot has only a certain capacity for the number of inspections that can be performed per day.

An optimal major maintenance plan is one that can supply adequate train-set for operation during periods of peak demand and limitations with minimum maintenance cost throughout the decision horizon. In CR system, experienced practitioners can generate effective and feasible plan. However, this process is time consuming, and an optimal solution is not guaranteed.

**Methodology**

In this section, a proposed approach and a corresponding optimization model are presented. As the optimization objective of EMUMMS is to ensure supply enough train-set to fulfill transportation tasks, not involved with the train-set operational plan, we divide the status of train-set into available status and major maintenance status (including the third-, fourth- and fifth-level maintenance status). The available status includes operation status, first- and second-level maintenance status. The reason not to set a waiting status is that a floating range existed which is helpful to flexibly adjust the start time of major maintenance.

**The Proposed Approach**

A space-time network is constructed based on graph theory, to analyze EMUMMS. In CR system, the running time of train-set is from 5:00 to 24:00 every day, and from 0:00 to 5:00 is for maintenance and recess accordingly. Furthermore, the working time of maintenance workshop is also in the duration of 5:00 to 24:00. Therefore, we construct a space-time network by setting the period from 0:00 to 5:00 as one node, and regarding the duration from 5:00 to 24:00 as arc that connects adjacent nodes (see Fig.1). The horizontal direction represents the elapse of time, and the vertical direction shows different status of train-set.

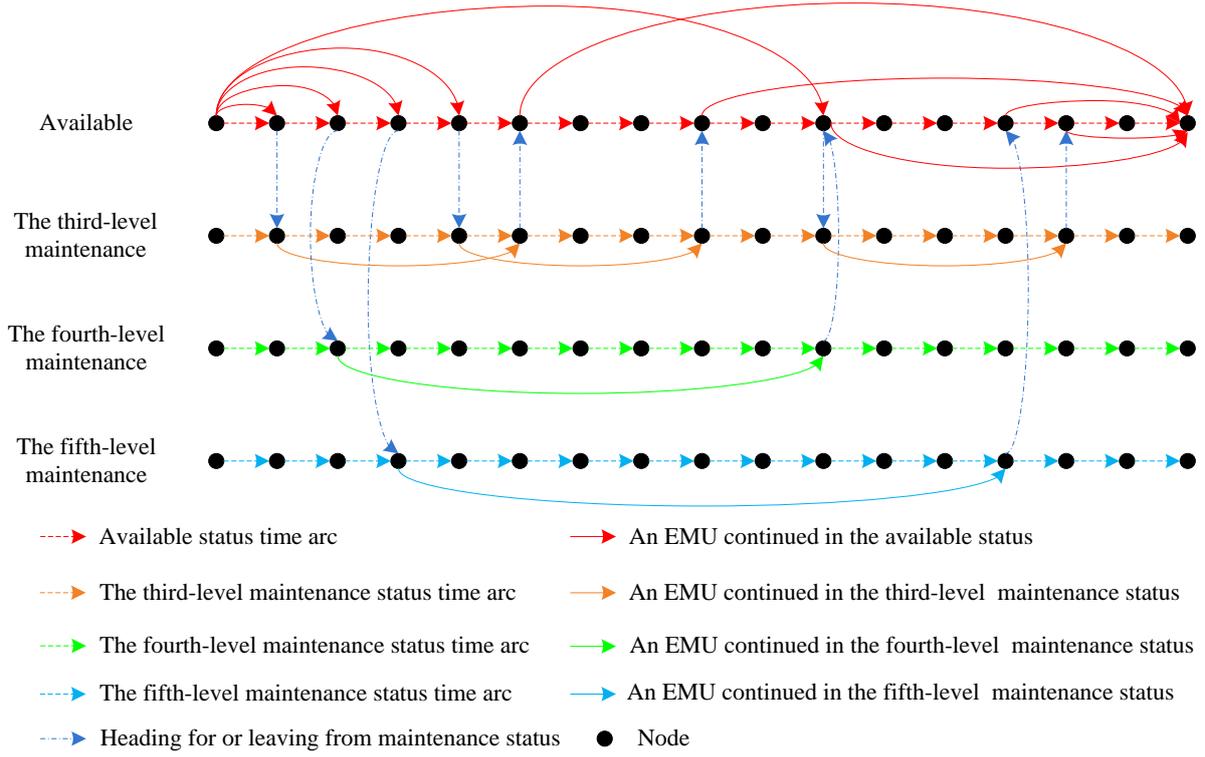

Fig. 2. Schematic diagram of connection relationship between available and major maintenance status for EMU

Let $V$ denote the set of all node and let $A$ denote the set of all arc. Directed graph $G$ consists of point set $V$ and arc set $A$, $G=(V,A)$. The last column denotes virtual super-node in planning horizon. The notations, including the indices, parameters, and sets, used above are listed in Table 3.

**Table 3.** Notation of indices, sets, and parameters.

| Notation | Description |
| --- | --- |
| *Indices* | |
| $i, j$ | The status of train-set; $i, j = 1, 2, 3, 4$, denote available status, the third-, fourth- and fifth-level maintenance status, respectively |
| $k, s$ | The ordinal number of column |
| $e$ | The train-set; $e \in E$ |
| *Parameters* | |
| $v_{ik}$ | The node locates in the $i$ row and the $k$ column; $v_{ik} \in V$ |
| $a$ | The arc; $a \in A$ |
| $a_{ik}^t$ | The time arc following $v_{ik}$; $a_{ik}^t \in A_i^t$ |
| $a_{ij}^k$ | The connecting arc from $a_{ik}^t$ to $a_{jk}^t$; $a_{ij}^k \in A^c$ |
| $a_{hk}^i$ | The sample arc presents the scope of train-set in the same status that from $a_{ih}^t$ to $a_{ik}^t$; $a_{hk}^i \in A^s$ |
| $w_{ik}$ | The weight of time arc |
| $n$ | The minimum value of $w_{1k}$ |
| $b_i$ | The minimum value of $w_{i+1,k}$, $i = 1, 2, 3$ |
| *Sets* | |
| $E$ | Set of train-set |
| $A$ | Set of all arc |
| $V$ | Set of all node |
| $A_i^t$ | Set of time arc;  $A_i^t \subset A^r$ |
| $A^c$ | Set of arc that connects different kinds of time arc;  $A^c \subset A^r$ |
| $A^r$ | Set of time arc and its connecting arc;  $A^r \subset A$ |
| $A^s$ | Set of sample arc;  $A^s \subset A$ |

We define the arc showed by dashed line as time arc in this paper. There are four kinds of time arcs exist in this network (distinguished by color), which represent available status, third-level, fourth-level and fifth-level maintenance status, respectively. Give these time arcs weights according to the practical meaning. Assign $w_{ik}$ to $a_{ik}^t$ that presents the number of train-sets on the $k$-th day and the $i$-th status. To successfully accomplish the transportation tasks of on the $k$-th day, the $w_{1k}$ must greater than $n$, i.e. $w_{1k} \geq n$. Note that the value of $n$ is different with train-set demands during different periods. $w_{2k}$ denotes the number of train-sets in third-level maintenance status on the $k$-th day. Due to the restricted service capacity of workshop, $w_{2k}$ has its maximum value $b_1$. Likewise, $w_{3k} \leq b_2$ and $w_{4k} \leq b_3$.

In this space-time network, to show a special status of train-sets intuitively, we connect the time arcs representing the same status to form an arc chain, depicted by a solid line painted with a same color. At the end of a planning period, direct all the solid line arcs to the super node at the end of each row and let $a_{hk}^i$ denote the solid arcs. In Fig.1, a train-set has been in available status from $v_{1,1}$ to $v_{1,5}$, which can be depicted by a red solid line from $v_{1,1}$ to $v_{1,5}$. The length of this red line denotes the duration of available status.

Similarly, this train-set has been in third-level maintenance status from $v_{2,5}$ to $v_{2,9}$, which can be denoted by an orange solid line from $v_{2,5}$ to $v_{2,9}$. The length of this line represents the time consumed of conducting third-level maintenance. It should be pointed out that Fig.1 is only a sketch map, whose time span is not true in practice. In addition, the process of train-set heading for and leaving from workshop are relatively stochastic in practical, so we do not take time consumption in heading and leaving process into consideration. For simplification, it is assumed that train-set can switch to major maintenance status on the same day after operation ending, and vice versa, which is depicted by deep-blue dash-dotted arc in Fig.1. The direction of arc signifies the status transition direction of train-set, denoted by $a_{ij}^k$, where $i$ and $j$ respectively present the status before and after transition. $k$ is the day of this process and let $A^c$ denote the set of these deep-blue dash dotted arcs. In this way, time arcs describing different status of a certain train-set are connected by deep-blue dash dotted arcs, forming a path throughout the planning cycle.

On the network, $a_{1k}^t$ and $a_{2k}^t$, $a_{1k}^t$ and $a_{3k}^t$, $a_{1k}^t$ and $a_{4k}^t$ can be linked together, respectively, by $a_{ij}^k$ ($i, j = 1, 2, 3, 4$). But it cannot be connected between the $a_{2k}^t$, $a_{3k}^t$ and $a_{4k}^t$. In other words, the status of train-set could transform from available to maintenance, or from maintenance to available, but cannot transform from a certain level of maintenance to another level.

EMUMMS determines which day the train-set should go to workshop and which level of maintenance the train-set should undertake. In a word, it is a problem of the status transformation between available and maintenance. On the network, this process can be shown as a path by connecting different $a$ ($a \in A^r$) from the beginning to the ending of planning horizon. Such as the path $r_0$: $v_{1,1} \xrightarrow{a_{1,1}^t} v_{1,2} \xrightarrow{a_{1,2}^t} v_{1,3} \xrightarrow{a_{1,3}^t} v_{1,4} \xrightarrow{a_{1,4}^t} v_{1,5} \xrightarrow{a_{1,2}^5} v_{2,5} \xrightarrow{a_{2,5}^t} v_{2,6} \xrightarrow{a_{2,6}^t} v_{2,7} \xrightarrow{a_{2,7}^t} v_{2,8} \xrightarrow{a_{2,8}^t} v_{2,9} \xrightarrow{a_{2,1}^9} v_{1,9} \xrightarrow{a_{1,9}^t} v_{1,10} \xrightarrow{a_{1,10}^t} v_{1,11} \xrightarrow{a_{1,11}^t} v_{1,12} \xrightarrow{a_{1,12}^t} v_{1,13} \xrightarrow{a_{1,13}^t} v_{1,14} \xrightarrow{a_{1,14}^t} v_{1,15} \xrightarrow{a_{1,15}^t} v_{1,16} \xrightarrow{a_{1,16}^t} v_{1,17}$, as shown in Fig. 1, describes that a certain train-set is in available status from the first day to the fourth day, then undertakes the third-level maintenance based maintenance requirements from the fifth day to the eighth day. After that, this train-set go back to undertake transport tasks from the ninth day to the end of planning horizon. In this way, all the status and status transition process of a certain train-set in a planning cycle are depicted by one path, and different status transition time point leads to different path. Therefore, convert EMUMMS problem into pathway chosen problem on the time-space network. At last, we optimize pathway chosen problem to realize the optimization of EMUMMS.

**Optimization model**

According to the proposed approach, an optimization model is constructed based on classic arc-path model. Let $r$ denote the feasible path for a certain train-set: $e$ ($e \in E$), and let $R_e$ denote the path set including all feasible paths on the time-space network. Based on the characteristics of arc-path model, we define the decision variables as follows:

$$x_e^r = \begin{cases} 1 & e \text{ chooses the } r-\text{path}, e \in E, r \in R_e \\ 0 & \text{otherwise} \end{cases}$$

And define the associated variables as follows:

$$\delta_e^{ra} = \begin{cases} 1 & a\text{-arc on the } r\text{-path of } e, e \in E, r \in R_e, a \in A^{ss} \\ 0 & \text{otherwise} \end{cases}$$

The associated variables describe inclusion relation between arc and path.

The objective function of the optimization model is to minimize the costs of major maintenance. We don't consider the optimization of certain major maintenance processes and regard it as a fixed value in this paper. So, the objective function of the optimization model can be solved by decreasing major maintenance frequency. And the practical method of decreasing major maintenance frequency is to increase the interval of two adjacent major maintenances, that is to say, to increase the actual operating mileage (and operating days) between two adjacent major maintenances. Moreover, according to maintenance requirements of each level, the actual operating mileage (and actual operating days) since the last major maintenance cannot exceed the upper limit of floating range. Therefore, the objective function of the optimization model could be described by minimizing the D-value between the cumulative operating mileage (and cumulative operating days) of maintenance requirements and the actual operating mileage (and operating days) since the last major maintenance.

According to the mentioned above, Let $L_e$ denote the value of cumulative operating mileage of major maintenance requirement for $e$ ($e \in E$). Let $\overline{l_e}$ denote the average daily mileage for the type of $e$. The average daily mileage of EMU train-set is statistic, and it is relatively precise to predict the eta of major maintenance because of the statistical time is long. The initial value $l_e^0$ and $t_e^0$ denote the accumulative operating mileage and operating days since last major maintenance. If a train-set is undergoing inspection at the beginning of planning horizon, $l_e^0$ and $t_e^0$ will take negative value for simplifying the model in this paper, and $t_e^0$ presents that $|t_e^0|$ days are remained to finish the inspection at this moment. And the $l_e^0$ equals to $t_e^0 \bar{l_e}$ accordingly. Let $l_e$ denote the status variable that indicates the expected operating mileage since last major maintenance when $e$ ($e \in E$) be sent to workshop. So, we have the following equation.

$$l_e = l_e^0 + (k-1)\sum_{r \in R_e} x_e^r \delta_e^{ra} \bar{l_e} \quad , \quad \forall e \in E , \tag{1}$$

where $a \in A^c$, i.e., $a_{1j}^k$ ($j = 2,3,4$). And $k$ is derived from $a_{1j}^k$ ($j = 2,3,4$). Such as $r_0$: $k = 5$.

In conclusion, the objective of the optimization model can be described as follows:

$$\min \quad Z = c\sum_e \max\{0, L_e - l_e\} \tag{2}$$

where $c$ denotes the unit costs of mileage loss resulting from the inspection ahead of schedule. Constraints of the model are mainly derived from the following several aspects.

(1) According to the proposed approach, any $e$ ($e \in E$) can discretionarily choose one and only one path from the space-time network, which is the uniqueness constraint of path in classical arc-path model. Its formulation is listed as follows:

$$\sum_{r \in R_e} x_e^r = 1 \quad , \quad \forall e \in E \tag{3}$$

(2) The minimum number of train-sets in available status is different based the demand of passengers, and the minimum number of train-sets in available status should be guaranteed. According to the value of $n$ during different periods, $a_{1k}^t$ can be divided into four stages including the usual, Spring Festival travel rush, summer holiday and National Day in this paper. Let $A_{1,1}^t$, $A_{1,2}^t$, $A_{1,3}^t$ and $A_{1,4}^t$ denote the subset of $A_1^t$ ($A_1^t = A_{1,1}^t + A_{1,2}^t + A_{1,3}^t + A_{1,4}^t$) on different stages and the minimum weights are $n_1$, $n_2$, $n_3$ and $n_4$, respectively. To describe this matter, a set of constraints established as follows:

$$\sum_{r \in R_e}\sum_{e \in E} x_e^r \delta_e^{ra} \geq n_1 , \quad \forall a \in A_{1,1}^t \tag{4}$$

$$\sum_{r \in R_e}\sum_{e \in E} x_e^r \delta_e^{ra} \geq n_2 , \quad \forall a \in A_{1,2}^t \tag{5}$$

$$\sum_{r \in R_e}\sum_{e \in E} x_e^r \delta_e^{ra} \geq n_3 , \quad \forall a \in A_{1,3}^t \tag{6}$$

$$\sum_{r \in R_e}\sum_{e \in E} x_e^r \delta_e^{ra} \geq n_4 , \quad \forall a \in A_{1,4}^t \tag{7}$$

(3) The number of train-sets in any level of maintenance status should not exceed the service capacity. To describe this matter, a set of constraints established as follows:

$$\sum_{r \in R_e}\sum_{e \in E} x_e^r \delta_e^{ra} \leq b_2 , \quad \forall a \in A_2^t \tag{8}$$

$$\sum_{r \in R_e} \sum_{e \in E} x_e^r \delta_e^{ra} \leq b_3 , \quad \forall a \in A_3^t \tag{9}$$

$$\sum_{r \in R_e} \sum_{e \in E} x_e^r \delta_e^{ra} \leq b_4 , \quad \forall a \in A_4^t \tag{10}$$

(4) Due to the floating range of cumulative operating mileage, presented in description of the problem, the train-set could be sent to workshop at flextime. But, it is obviously unreasonable to be too early or too late. Let $\xi_e^l$ denote the lower bound of floating range and $\xi_e^u$ denote the upper bound for $e$ ($e \in E$). This matter can be described by inequality:

$$\xi_e^l \leq l_e \leq \xi_e^u, \quad \forall e \in E \tag{11}$$

From the above analysis, the model for optimizing EMUMMS is summarized as model I.

$$\min \quad Z = c \sum_e \max\{0, L_e - l_e\}$$

$$s.t. \quad \text{Constraints (3)} \sim \text{(11)}$$

$$x_d^r \in \{0, 1\}$$

Moreover, let $l_e'$ denote the total operating mileage for $e$ ($e \in E$). And the level of last maintenance is denoted by $M'$.

**Algorithm design**
**The generation of feasible path set**

To solve this model, we have to generate $R_e$ for every $e$ ($e \in E$) on the network firstly. According to the network design, the path is an arc chain that formed by a variety end-to-end arcs throughout the whole planning horizon. To find out these arc chains, we should recognize the subsequent arc set of each type of arc. At the beginning of a planning horizon, the train-set may be in available or major maintenance status.

The subsequent arc set of $a_{1k}^t$ is described as follows: As for the train-set is available, the subsequent status of it may be still available or turn to major maintenance. So the subsequent arc set of $a_{1k}^t$ is { $a_{1,k+1}^t$, $a_{1j}^{k+1}$ ($j=2,3,4$) }.

To understand the definition of the subsequent arc set of $a_{1j}^k$ ($j=2,3,4$), we can consider this case: along with the train-set go to workshop, the maintenance procedures of a certain level will be conducted on it. Thus, the subsequent arc set of $a_{1j}^k$ ($j=2,3,4$) is { $a_{jk}^t$ ($j=2,3,4$) }.

The subsequent arc set of $a_{jk}^t$ ($j=2,3,4$) is illustrated as follows: when a certain level of major maintenance for train-set is completed, the status of train-set will change to available; but, if further maintenance is needed, the subsequent arc should be $a_{j,k+1}^t$ ($j=2,3,4$). Therefore, { $a_{j1}^{k+1}$ ($j=2,3,4$), $a_{j,k+1}^t$ ($j=2,3,4$) } is the subsequent arc set of $a_{jk}^t$ ($j=2,3,4$).

The subsequent arc set of $a_{j1}^k$ ($j=2,3,4$) is defined as follows: once the major maintenance of train-set is completed, the status of train-set will change to available. So the subsequent arc set of $a_{j1}^k$ ($j=2,3,4$) is the { $a_{1k}^t$ }.

After clearing the subsequent arc sets of all arcs on the network, $R_e$ can be generated.

Moreover, $R_e$ can be downsized according to the constraint (11). So, a feasible path set denoted by $R'_e$ can be generate by depth-first search algorithm, the specific steps are as follows:

Step1: According to the initial status of train-set, we assign the initial arc and clear the arc's type.

Step2: Take constraint (11), $l_e$, $M'$ and maintenance requirements as criterion. Based on the connection sequence of arcs mentioned above, the depth-first search algorithm is applied to find the subsequent arc and clear the type of subsequent arc. If none subsequent arc is found, then turn to Step 5. Otherwise, judge whether or not the subsequent arc points at the virtual super-node. If yes, turn to Step 3; otherwise, turn to Step 2.

Step 3: Connect all arcs in the order and output the path: $r$. If $r$ already exists, turn to Step5; otherwise, turn to Step 4.

Step 4: Add $r$ to $R'_e$ ( $R'_e \subset R_e$ ), and delete the $a_{1j}^k$ ( $j=2,3,4$ ) and $a_{j1}^k$ ( $j=2,3,4$ ) exist in path $r$ from the network. Then reconstruct the network $G'$ and go back to Step1.

Step5: End the feasible path generation process and output $R'_e$.

Reapply the algorithm above and generate $R'_e$ on the network for every $e$ ( $e \in E$ ). And, let $R$ denote the set of all feasible path, so, $R = \sum_{e \in E} R'_e$.

As constraint (11) has already been used in the process of generating $R$, the optimization model of major maintenance schedule can be simplified as model II:

$$\min \quad Z = c \sum_e \max\{0, L_e - l_e\}$$

$$\text{s.t.} \quad \sum_{r \in R'_e} x_e^r = 1, \quad \forall e \in E \tag{13}$$

$$\sum_{r \in R'_e} \sum_{e \in E} x_e^r \delta_e^{ra} \geq n_1, \quad \forall a \in A_{1,1}^t \tag{14}$$

$$\sum_{r \in R'_e} \sum_{e \in E} x_e^r \delta_e^{ra} \geq n_2, \quad \forall a \in A_{1,2}^t \tag{15}$$

$$\sum_{r \in R'_e} \sum_{e \in E} x_e^r \delta_e^{ra} \geq n_3, \quad \forall a \in A_{1,3}^t \tag{16}$$

$$\sum_{r \in R'_e} \sum_{e \in E} x_e^r \delta_e^{ra} \geq n_4, \quad \forall a \in A_{1,4}^t \tag{17}$$

$$\sum_{r \in R'_e} \sum_{e \in E} x_e^r \delta_e^{ra} \leq b_2, \quad \forall a \in A_2^t \tag{18}$$

$$\sum_{r \in R'_e} \sum_{e \in E} x_e^r \delta_e^{ra} \leq b_3, \quad \forall a \in A_3^t \tag{19}$$

$$\sum_{r \in R'_e} \sum_{e \in E} x_e^r \delta_e^{ra} \leq b_4, \quad \forall a \in A_4^t \tag{20}$$

$$x_d^r \in \{0, 1\} \tag{21}$$

**Solution strategy of the model**

Genetic algorithm is good at controlling the search process but may easily result in premature convergence for the problem. While the simulated annealing algorithm performs well in local search but performs bad in global search. Therefore, integrate these two algorithms and then design a new solution approach according to the characteristics of our model.

(1) Constraint removing

As for constraint (3) in model II, we can realize it by setting coding principle for chromosomes, which will be described in detail. Because the decision variables are 0-1 variables, the binary coding approach is proposed. Each chromosome represents a solution of the model, i.e., the major maintenance schedule of all train-sets. The length of chromosome is $|R| = \sum |R'_e|$ ($|R'_e|$ denotes the number of feasible paths of $e$). And the block coding approach is applied to chromosomes. Use "|" to divide the chromosome into blocks that presents a feasible path set for certain train-set. The number of blocks equal to the number of train-sets. (see Fig. 3).

$$\left| x_1^1 x_1^2 \cdots x_1^{|R_1|} \right| x_2^1 x_2^2 \cdots x_2^{|R_2|} \left| \cdots \right| x_e^1 x_e^2 \cdots x_e^{|R'_e|} \right|$$

Fig. 3. Diagram of chromosome encoding

The gene sequence in each block corresponding to the order of paths in the feasible path set. If train-set chooses a path, the corresponding gene will be set to 1. Moreover, according to constraint (3), the total value of genes in each block equals 1. That is to say, only one gene in a block can equal1, and the rest can only be 0.

Through the penalty function and coding approach mentioned above, the model II can be converted to an optimization problem with no constraints.

$$\min \quad Z = c \sum_e \max\{0, L_e - l_e\} + \lambda_1 \sum_{i=1}^{4} \sum_{a \in A_{1,i}^t} \max\{0, n_i - \sum_{r \in R'_e} \sum_{e \in E} x_e^r \delta_e^{ra}\} + \lambda_2 \sum_{i=2}^{4} \sum_{a \in A_i^t} \max\{0, \sum_{r \in R'_e} \sum_{e \in E} x_e^r \delta_e^{ra} - b_i\} \quad (22)$$

where $c$ is $0.0001$, $\lambda_1 = 100$ and $\lambda_2 = 80$. Since the constraint of insufficient train-set is more intense than that of maintenance service capacity, the value of $\lambda_1$ is greater than $\lambda_2$.

(2) Relevant strategies for genetic algorithm with simulated annealing survival mechanism

a. Strategy for chromosome coding and initial solution generation

The chromosome coding approach has been previously described in detail, so this section will focus on the strategy for initial solution generation. Penalty function and appropriate coding approach are applied to model II and convert the problem into unconstrained problem. Therefore, the generation of initial solution only needs to meet the requirements of the unconstrained problem. According to the coding approach, set certain gene to 1 and set others to 0 in each coding block. In this way, an initial solution can be obtained. Then, repeat this procedure and generate the initial population with required size, such as $sizepop = 30$. It is inevitable that some solutions in initial population do not meet the constraints of model I. Despite this, due to the strong convergence and robust of genetic algorithm, the generation of the optimization solutions through continuous evolution will not be affected.

b. Calculation strategy for individual fitness

Fitness calculation function is developed to evaluate the individual fitness in iteration.

$$Fit_j^t = 1/Z_j^t \tag{23}$$

In equation (23), $Fit_j^t$ denotes the fitness of chromosome $j$ after $t$ times of iteration; $Z_j^t$ is the objective function value of chromosome $j$ after $t$ times of iteration. In this way, the better quality of chromosome, the greater fitness, i.e., the smaller value of objective function, which accords with the laws of evaluation.

c. Selection strategy

The fitness function is taken as the evaluation criteria in the selection of chromosomes. Therefore, the chromosome has a high probability be selected with the great value of fitness. The roulette wheel selection strategy in standard genetic algorithm is adopted to select chromosomes. The probability of an individual to be selected is as follows:

$$p_j = Fit_j^t / \sum_j Fit_j^t \tag{24}$$

The roulette wheel is divided into different regions according to the probability we got. Then, a random number ($\zeta_1$) is generated in [0, 1] and chromosomes are selected based on the region where $\zeta_1$ appears. Select repeatedly until the number of chromosomes reaches $sizepop$, the size of population, based on the roulette wheel selection strategy.

d. Crossover strategy

Let $p_c$ denote the crossover probability and set it as 0.88 in this paper. Pair the chromosomes randomly, and adopt the genetic crossover approach in units of blocks. The detailed strategy is as follows. A random number ($\zeta_2$) is generated firstly. Then, if $\zeta_2 < p_c$, randomly produce two integers in [1, $|R|$] ($\omega_1$ and $\omega_2$, $\omega_1 < \omega_2$) and swap the gene blocks from $\omega_1$ to $\omega_2$ between two parent chromosomes in a pair. In this way, two child chromosomes be obtained.

e. Mutation strategy

Let $p_m$ denote the mutation probability and set it as 0.08 in this paper. Just like crossover operation, the mutation operation for chromosome is conducted in units of blocks, too. Firstly, generate a random number ($\zeta_3$) for a certain chromosome. If $\zeta_3 < p_m$, generate a random integer ($\omega_3$) in interval [1, $|R|$]. Secondly, generate another random integer ($\omega_4$) in interval [1, $|R_{\omega_3}|$]. And if the value of gene $\omega_4$ in gene block $\omega_3$ equals to 1, regenerate $\omega_4$ until the value of gene $\omega_4$ equals to 0. Finally, set the value of gene $\omega_4$ to 1 and set the rest genes in block $\omega_3$ to 0.

f. Survival strategy of simulated annealing

In genetic algorithm, the parent chromosomes will be replaced by the offspring after selection, crossover and mutation. However, regular updating methods are likely to produce local optimum solution. To overcome this disadvantage, simulated annealing survival strategy should be adopted for updating. The detailed procedures are as follows. If the fitness of an offspring individual is better than its parents, accept this offspring individual undoubtedly. Otherwise, accept the offspring individual with probability: $\exp\left(\left(Fit_j^{t+1} - Fit_j^t\right)/T^l\right)$. Where $T^l$ denotes the temperature after $l$ times of cooling in simulated annealing process. $T^{l+1} = \alpha T^l$ and $\alpha$ is the cooling rate is selected as 0.8 in this paper.

**Algorithm procedure**

The flowchart of the algorithm is shown in Fig. 4.

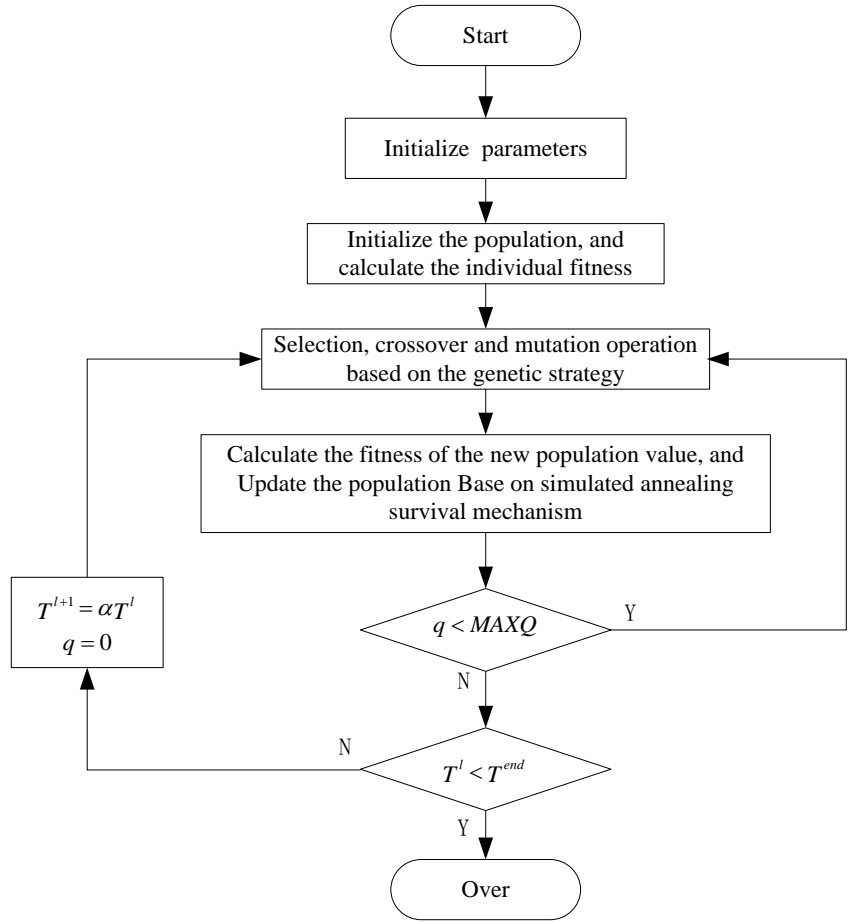

Fig. 4. The flowchart of the algorithm

The specific steps are as follows:

Step1: Initialize the maximum evolution times $MAXT = 3000$, original temperature $T^0 = 100$, ending temperature $T^{end} = 1$. and set the iteration time as $t = 0$.

Step 2: Generate original population randomly based on coding strategy and compute the fitness of individuals $Fit_j^t$ $j = 1, 2, \cdots, sizepop$.

Step3: Generate a new population through selection, crossover and mutation according to relevant genetic strategies.

Step 4: Compute the fitness of individuals in the new population $Fit_j^{t+1}$, $j = 1, 2, \cdots, sizepop$. If $Fit_j^{t+1} > Fit_j^t$, replace the parents with offspring individuals; otherwise, accept offspring individuals with probability $\exp\left((Fit_j^{t+1} - Fit_j^t)/T^l\right)$, $t = t + 1$.

Step5: If $t < MAXT$, go back to Step3; otherwise, turn to Step 6.

Step6: If $T^l < T^{end}$, end the algorithm and output the optimum solution; otherwise, set

$T^{l+1} = \alpha T^l$, $t = 0$ and operate the simulated annealing procedures, then go back to Step 3.

**Conclusion**

This research optimizes the electric multiple unit major maintenance schedule from the perspective of network design by transforming it into the pathway chosen problem on the time-space network. A 0-1 integer programming model is developed considering the different transportation demand for passengers and the capacity limited of workshop, etc. Compared with the manual process, the genetic algorithm with simulated annealing survival mechanism is also developed to improve solution quality and efficiency. Using this decision support tool can help railways with similar characteristics to improve the efficiency in electric multiple unit major maintenance schedule. What have been done in this research is to optimize electric multiple unit major maintenance strategically. Future research may investigate the possibility to optimize the maintenance plan for the first-, second-level inspection considering the different transportation demand for passengers and the maintenance capacity limited of sub-depot under the finite fleet-size of train-set.

**References**

1. https://en.wikipedia.org/wiki/China_Railway_High-speed.
2. Stuchly, V., Grencik, J., Poprocky, R.,. Railway vehicle maintenance and information systems. *Computers in railways VII: International conference on computers in railways*, 2000, pp.885-894. Bologne, the France.
3. Shimada N. Rolling stock maintenance for safe and stable transport. *Japanese Railway Eng*, 2006; 46(2): 4-7.
4. Cheng, Y. H., Tsao, H. L. Rolling stock maintenance strategy selection, spares parts' estimation, and replacements' interval calculation. *Int J Prod Econ*, 2010; 128(1): 404-412.
5. Sriskandarajah, C., Jardine, A. K. S., Chan, C. K. Maintenance scheduling of rolling stock using a genetic algorithm. *J Oper Res Soc* 1998; 49(11): 1130–1145.
6. Maróti, G., Kroon, L. Maintenance routing for train units: the transition model. *Transp Sci* 2005; 39(4): 518-525.
7. Maróti, G., Kroon, L. Maintenance routing for train units: the interchange model. *Comput Oper Res* 2007; 34(4): 1121-1140.
8. Alfieri, A., Groot, R., Kroon, L. G., et al. Efficient circulation of railway rolling stock. *Transp Sci* 2006; 40(3): 378–391.
9. Rezvanizaniani1, S. M., Valibeigloo1, M., Asghari1 M, et al. Reliability centered maintenance for rolling stock: a case study in coaches' wheel sets of passenger trains of Iranian railway. *2008 IEEE International Conference on Industrial Engineering and Engineering Management*, 2008, pp.516-520. Washington.
10. Tsuji, Y., Kuroda, M., Imoto., Y. Rolling stock planning for passenger trains based on ant colony optimization. *Transactions of the Japan Society of Mechanical Engineers C* 2010; 76: 397-406.
11. WANG, Y., LIU, J., MIAO, J. R. Column generation algorithms based optimization method for

maintenance scheduling of multiple units. *China Railway Sci* 2010; 31(2): 115-120. (in Chinese).

12. Giacco, G. L., D'Ariano, A. Pacciarelli, D. Rolling stock rostering optimization under maintenance constraints. *J Int Transp Sys* 2014; 18(1): 95–105.
13. Lai, Y. C., Fan, D. C., Huang, K. L. Optimizing rolling stock assignment and maintenance plan for passenger railway operations. *Comput ind eng* 2015; 85: 284–295.
14. Ziaratia, K., Soumisa, F., Desrosiers, J., et al. Locomotive assignment with heterogeneous consists at CN North America. *Eur J Oper Res* 1997; 97(2): 281-292.
15. Lingaya, N., Cordeau, J. F., Desaulniers, G., et al. Operational car assignment at VIA Rail Canada. *Transp Res Part B: Methodol* 2002; 36(9): 755–778.
16. WANG, L., MA, J. J., LIN, B. L., et al. Method for optimization of freight locomotive scheduling and routing problem. *J China Railway Soc* 2014; 36(11): 7-15.
17. Moudania, W. E., Félix, M. C. A dynamic approach for aircraft assignment and maintenance scheduling by airlines. *J Air Transp Manage* 2000; 6(4): 233-237.
18. Budai, G., Huisman, D., Dekker, R. Scheduling preventive railway maintenance activities. *J Oper Res Soc* 2006; 57(9): 1035–1044.
19. Mehmet, B., Bilge, Ü. Operational aircraft maintenance routing problem with remaining time consideration. *Eur J Oper Res* 2014; 235(1): 315-328.
20. Grigoriev, A., Klundert, J. V. D., Spieksma, F. C. R. Modeling and solving the periodic maintenance problem. *Eur J Oper Res* 2006; 172(3): 783–797.
21. China Railway. User Manual for EMU Operation and Maintenance Procedures 2013.